\documentclass[10pt,reqno]{amsart}

\usepackage{a4wide}
\usepackage{dsfont}                   
\usepackage{amssymb}
\usepackage[all]{xy}
\usepackage{mathrsfs}                 

\newtheorem{proposition}{Proposition}[section]
\newtheorem{theorem}[proposition]{Theorem}

\theoremstyle{remark}
\newtheorem{definition}[proposition]{Definition}
\newtheorem{example}[proposition]{Example}

\newcommand{\cst}{\ifmmode\mathrm{C}^*\else{$\mathrm{C}^*$}\fi}
\newcommand{\st}{\;\vline\;}
\newcommand{\CC}{\mathbb{C}}
\newcommand{\ZZ}{\mathbb{Z}}
\newcommand{\vt}{\!\vartriangle\!}
\newcommand{\XX}{\mathbb{X}}
\newcommand{\YY}{\mathbb{Y}}
\newcommand{\EE}{\mathbb{E}}
\newcommand{\DD}{\mathbb{D}}
\newcommand{\KK}{\mathbb{K}}
\newcommand{\sS}{\mathbb{S}}
\newcommand{\HH}{\mathbb{H}}
\newcommand{\GG}{\mathbb{G}}
\newcommand{\WW}{\mathbb{W}}
\newcommand{\PP}{\mathbb{P}}
\newcommand{\FF}{\mathbb{F}}
\newcommand{\UU}{\mathbb{U}}
\newcommand{\MM}{\mathbb{M}}
\newcommand{\sA}{\mathsf{A}}
\newcommand{\sB}{\mathsf{B}}
\newcommand{\sC}{\mathsf{C}}
\newcommand{\cH}{\mathscr{H}}
\newcommand{\tens}{\otimes}
\newcommand{\id}{\mathrm{id}}
\newcommand{\comp}{\!\circ\!}
\renewcommand{\Bar}[1]{\overline{#1}}
\newcommand{\I}{\mathds{1}}
\newcommand{\Label}[1]{\label{#1}\stepcounter{equation}\tag{\theequation}}
\DeclareMathOperator{\C}{C}
\DeclareMathOperator{\B}{B}
\DeclareMathOperator{\M}{M}
\DeclareMathOperator{\Mor}{Mor}
\DeclareMathOperator{\cK}{\mathcal{K}}

\numberwithin{equation}{section}

\begin{document}

\author{Piotr M.~So{\l}tan} \address{Department of Mathematical Methods in Physics, Faculty of Physics, University of Warsaw, Poland}
\email{piotr.soltan@fuw.edu.pl}

\title{Quantum semigroup structures on quantum families of maps}

\thanks{Partially supported by National Science Centre (NCN) grant no.~2011/01/B/ST1/05011.}



\begin{abstract}
These are notes accompanying lectures at 7th ILJU School of Mathematics ``Banach Spaces and Related Topics" in Gyeongju, Koera. The lectures are devoted to exposition of the theory of quantum families of maps with emphasis on the study of quantum semigroup structures appearing in this context.
\end{abstract}

\date{January 28 -- February 1, 2013}

\maketitle


\section*{{\Large Lecture I}}
\renewcommand{\thesection}{I.\arabic{section}}

\section{C*-algebras and Gelfand's duality}

A \emph{Banach algebra} is a Banach space $\sA$ endowed with a bilinear operation
\[\Label{prod}
\sA\times\sA\ni(a,b)\longmapsto{ab}\in\sA
\]
such that $(ab)c=a(bc)$ for all $a,b,c\in\sA$ and
\[
\|ab\|\leq\|a\|\|b\|
\]
for all $a,b\in\sA$. The bilinear operation \eqref{prod} is called a \emph{product} or \emph{multiplication} of $\sA$. A typical example of a Banach algebra is the space $\B(E)$ of all bounded linear operators on a Banach space $E$ with composition as product and the operator norm, but there are many other examples (see below). We will only consider Banach algebras over the field of complex numbers. In this context an \emph{involution} or \emph{$*$-operation} on a Banach algebra $\sA$ is an antilinear mapping
\[\Label{cstid}
\sA\ni{a}\longmapsto{a^*}\in\sA
\]
which satisfies
\begin{itemize}
\item $(a^*)^*=a$ for all $a\in\sA$ (involutivity),
\item $(ab)^*=b^*a^*$ for all $a,b\in\sA$ (anti-multiplicativity).
\end{itemize}
One usually adds a condition that $\|a\|=\|a^*\|$, i.e.~that the $*$-operation is isometric. A Banach algebra with an involution is called an \emph{involutive Banach algebra} or a \emph{Banach $*$-algebra}.

\begin{definition}
A \cst-algebra is a Banach $*$-algebra $\sA$ such that
\[
\|a^*a\|=\|a\|^2
\]
for any $a\in{A}$.
\end{definition}

Here are the most important examples of \cst-algebras:

\begin{example}\label{exCX}
Let $X$ be a compact space (by definition all compact and locally compact spaces are Hausdorff). Then the Banach space $\C(X)$ (with $\sup$ norm) is a \cst-algebra with multiplication defined as the pointwise multiplication of functions and involution $f\mapsto{f^*}$ given by
\[
f^*(x)=\Bar{f(x)},\qquad(x\in{X}).
\]
The \cst-algebra $\C(X)$ has a \emph{unit element}, namely the function $\I\colon{X}\ni{x}\mapsto{1}\in\CC$. More generally, if $Y$ is a locally compact space then the Banach space
\[
\C_0(Y)=\bigl\{f\in\C(Y)\st\text{ for any $\delta>0$ the set }\{y\in{Y}\st|f(y)|>\delta\}\text{ is compact}\bigr\}
\]
is a \cst-algebra with the same norm and algebraic structure as described for $\C(X)$. However, $\C_0(Y)$ does not have a unit unless $Y$ is compact (ant then $\C_0(Y)=\C(Y)$). The algebras $\C(X)$ and $\C_0(Y)$ are \emph{commutative}, i.e.~$fg=gf$ for all $f,g$ in $\C(X)$ (resp.~$\C_0(Y)$).
\end{example}

\begin{example}\label{exKH}
Let $\cH$ be a Hilbert space and consider the Banach space $\cK(\cH)$ of all \emph{compact} operators acting on $\cH$. Then with composition of operators as multiplication and Hermitian conjugation as the involution $\cK(\cH)$ is a \cst-algebra. It is unital (i.e.~has a unit) if and only if $\cH$ is finite-dimensional. Unless $\dim{\cH}=1$, $\cK(\cH)$ is not commutative.
\end{example}

Let $\sA$ and $\sB$ be \cst-algebras. A map $\Phi$ form $\sA$ to $\sB$ is called a \emph{$*$-homomorphism} if $\Phi$ is a homomorphism of algebras $\sA\to\sB$ and $\Phi(a^*)=\Phi(a)^*$ for all $a\in\sA$. It turns out that any $*$-homomorphism is automatically a contraction for the norms of $\sA$ and $\sB$. A $*$-isomorphism is a bijective $*$-homomorphism (which is automatically isometric).

It turns out that the \cst-algebras presented in Examples \ref{exCX} and \ref{exKH} are to some extent general. In the context of Example \ref{exCX} we have the famous theorem of Gelfand which will be discussed in detail below.

\begin{theorem}[I.M.~Gelfand]\label{IMG}
Let $A$ be a commutative \cst-algebra with unit. Then there exist a compact space $X$ and a $*$-isomorphism of $\C(X)$ onto $\sA$. The space $X$ is unique up to homeomorphism.
\end{theorem}

By a simple trick of \emph{adjoining the unit} one can extend Theorem \ref{IMG} to the case of non-unital \cst-algebras. More importantly, morphisms between \cst-algebras can be characterized in a similar way. In order to be more precise we will make a brief digression and define morphisms of \cst-algebras.

\subsection{Multipliers and morphisms}

Let $\sA$ be a \cst-algebra. The \emph{multiplier algebra} of $\sA$ is a canonically defined unital \cst-algebra $\M(\sA)$ such that $\sA$ is embedded as an ideal in $\M(\sA)$ which is \emph{essential}, i.e.~any other non-zero ideal of $\M(\sA)$ has an non-zero intersection with $\sA$. The defining property of $\M(\sA)$ is that if $\sC$ is a unital \cst-algebra and we are given an embedding of $\sA$ into $\sC$ as an essential ideal, then there exists an injective $*$-homomorphism of $\sC$ into $\M(\sA)$ extending the identity map $\sA\to\sA$ (considered as a map from a subset of $\sC$ onto a subset of $\M(\sA)$). It can be shown that $\M(\sA)=\sA$ if and only if $\sA$ is unital.

\begin{example}
Let $Y$ be a locally compact space. Then the multiplier algebra $\M\bigl(\C_0(Y)\bigr)$ of $\C_0(Y)$ is canonically isomorphic to the \cst-algebra $\C_{\text{b}}(Y)$ of all bounded continuous functions on $Y$. This algebra can be further shown to be isomorphic to the algebra of all continuous functions on the \emph{Stone-\v{C}ech compactification} of $Y$.
\end{example}

\begin{example}
Let $\cH$ be a Hilbert space. Then the multiplier algebra $\M\bigl(\cK(\cH)\bigr)$ of $\cK(\cH)$ is canonically isomorphic to the \cst-algebra $\B(\cH)$ of all bounded operators on $\cH$.
\end{example}

Now let $\sA$ and $\sB$ be \cst-algebras. A \emph{morphism} from $\sA$ to $\sB$ is a $*$-homomorphism $\Phi\colon\sA\to\M(\sB)$ such that $\Phi(\sA)\sB$ is dense in $\sB$. The notation $\Phi(\sA)\sB$ stands for the linear span of all products of the form $\Phi(a)b$ with $a\in\sA$ and $b\in\sB$. Recall that $\sB$ is an ideal in $\M(\sB)$ so such products are automatically elements of $\sB$. The set of all morphisms from $\sA$ to $\sB$ will be denoted by $\Mor(\sA,\sB)$.

Any $\Phi\in\Mor(\sA,\sB)$ admits a unique extension to a $*$-homomorphism $\Bar{\Phi}$ of unital algebras $\M(\sA)\to\M(\sB)$. It follows that if $\sC$ is another \cst-algebra and $\Psi\in\Mor(\sC,\sA)$ then we have a well defined map $\Bar{\Phi}\comp\Psi\colon\sC\to\M(\sB)$. It turns out that this is a morphism form $\sC$ to $\sB$. In particular morphisms can be composed.

With help of the notion of morphisms of \cst-algebras we can re-state the Gelfand theorem in its fuller version:

\begin{theorem}
\noindent
\begin{enumerate}
\item Let $\sA$ be a commutative \cst-algebra. Then there exists a unique up to homeomorphism locally compact space $X$ such that $\sA$ is isomorphic to the \cst-algebra $\C_0(X)$.
\item Let $\sA$ and $\sB$ be commutative \cst-algebras and let $X$ and $Y$ be locally compact spaces such that $\sB=\C_0(X)$ and $\sA=\C_0(Y)$. Then
\begin{enumerate}
\item any continuous $\phi\colon{X}\to{Y}$ defines a morphism $\Phi\in\Mor(\sA,\sB)$ via
\begin{equation}\label{Phif}
\Phi(f)=f\comp\phi,\qquad(f\in\sA);
\end{equation}
\item for any $\Phi\in\Mor(\sA,\sB)$ there exists a continuous $\phi\colon{X}\to{Y}$ such that \eqref{Phif} holds.
\end{enumerate}
\end{enumerate}
\end{theorem}

Thus any statement about locally compact spaces and continuous maps between them can be formulated as a statement about commutative \cst-algebras and their morphisms.

Before continuing let us note the following: let $\sA$ and $\sB$ be \cst-algebras. Among many norms on the algebraic tensor product $\sA\tens_{\text{alg}}\sB$ there is always the \emph{smallest} norm satisfying (the appropriate analog of) \eqref{cstid}. This norm is usually called the \emph{minimal} \cst-norm on $\sA\tens_{\text{alg}}\sB$ and the completion of $\sA\tens_{\text{alg}}\sB$ with respect to this norm, denoted by $\sA\tens\sB$ is called the \emph{minimal tensor product} of $\sA$ and $\sB$. In case either $\sA$ or $\sB$ is commutative (and also in many other cases) this is, in fact, the only \cst-norm on $\sA\tens_{\text{alg}}\sB$. Moreover we have

\begin{proposition}\label{cartprod}
Let $X$ and $Y$ be locally compact spaces. Then the minimal tensor product $\C_0(X)\tens\C_0(Y)$ is canonically isomorphic to the \cst-algebra $\C_0(X\times{Y})$. Under this isomorphism the tensor product $f\tens{g}$ with $f\in\C_0(X)$ and $g\in\C_0(Y)$ is mapped to the function
\[
X\times{Y}\ni(x,y)\longmapsto{f(x)g(y)}\in\CC.
\]
\end{proposition}

\section{Quantum spaces and their morphisms}

The Gelfand theorem tells us that the category of locally compact spaces and continuous maps is dual to the category of commutative \cst-algebras with morphisms. But the class of all (not necessarily commutative) \cst-algebras with morphisms is also a category (with the commutative ones forming a full subcategory). Can we treat all \cst-algebras as algebras of functions? Of course, a non-commutative \cst-algebra cannot be isomorphic to an algebra of functions on a space, but it can be \emph{regarded as one}. Thus we arrive at the following definition:

\begin{definition}
A \emph{quantum space} is an object of the category dual to the category of all \cst-algebras with morphisms.
\end{definition}

Thus any \cst-algebra $\sA$ is thought of as the algebra of functions on a \emph{quantum space} $\XX$. From now on we will write $\C_0(\XX)$ instead of $\sA$ and treat the virtual object $\XX$ in a way reminiscent of studying a locally compact space. We will call a quantum space $\XX$ \emph{compact} if the \cst-algebra $\C_0(\XX)$ is unital and in this case we will use the symbol $\C(\XX)$ to denote this \cst-algebra. Following this analogy we will say that a quantum space $\MM$ is \emph{finite} if $\C_0(\MM)$ is a finite-dimensional algebra (in this case it is automatically unital, so according to the remarks above we will denote it by $\C(\MM)$).

The concept of a map between quantum spaces is also quite clear. If $\XX$ and $\YY$ are quantum spaces (i.e.~we are given two \cst-algebras which we denote $\C_0(\XX)$ and $\C_0(\YY)$) then a continuous map from $\XX$ to $\YY$ is, by definition, an element of $\Mor\bigl(\C_0(\YY),\C_0(\XX)\bigr)$.

It has to be stressed that we are not developing a theory of some new kind of mathematical object. Any theorem about quantum spaces and their maps is a theorem about \cst-algebras and their morphisms. We end this section with the following warning: motivated by Proposition \ref{cartprod} one is tempted to define the Cartesian product of two quantum spaces $\XX$ and $\YY$ as the quantum space $\EE$ such that $\C_0(\XX)\tens\C_0(\YY)=\C_0(\EE)$. In fact we will almost do this. More precisely we will discuss various objects and phenomena in the category of quantum spaces in which the tensor product of \cst-algebras will play a role analogous to the role played by the Cartesian product of topological spaces. However this should not be taken as a strict definition because the categorical properties of the tensor product are not quite the same in the category of commutative algebras (where the tensor product is a categorical sum) and in the category of all \cst-algebras.

\section{Compact quantum groups and quantum semigroups}

\begin{definition}[S.L.~Woronowicz]\label{defCQG}
A \emph{compact quantum group} is a compact quantum space $\GG$ such that there exists a $\Delta\in\Mor\bigl(\C(\GG),\C(\GG)\tens\C(\GG)\bigr)$ such that
\begin{enumerate}
\item $(\Delta\tens\id)\comp\Delta=(\id\tens\Delta)\comp\Delta$ (coassociativity)
\item\label{defCQG2} the sets $\Delta\bigl(\C(\GG)\bigr)\bigl(\I\tens\C(\GG)\bigr)$ and $\bigl(\C(\GG)\tens\I\bigr)\Delta\bigl(\C(\GG)\bigr)$ are dense in $\C(\GG)\tens\C(\GG)$ (density conditions).
\end{enumerate}
\end{definition}

The map $\Delta$ form Definition \ref{defCQG} is called the \emph{comultiplication} or \emph{coproduct} while the density conditions in point
\eqref{defCQG2} are referred to as the \emph{cancellation properties}. The reason for this is that if $\GG$ (as in Definition \ref{defCQG}) is a classical space (i.e.~the \cst-algebra $\C(\GG)$ is commutative) then writing $G$ instead of $\GG$ we obtain form $\Delta$ a continuous associative two-argument operation $\mu\colon{G}\times{G}\to{G}$ on the compact space $G$ and the density conditions are equivalent to the fact that $\mu$ obeys cancellation from the right and from the left. It is a known fact that a compact semigroup whose multiplication obeys cancellation properties is a compact group, and so any compact group provides an example of a compact quantum group.

Compact quantum groups have been studied for many years now and many interesting examples were discovered. However, our main focus in these notes will be on more general objects. An example of such a more general structure is that of a compact quantum semigroup: a compact quantum space $\sS$ is a \emph{compact quantum semigroup} if it is equipped with a coassociative $\Delta\in\Mor\bigl(\C(\sS),\C(\sS)\tens\C(\sS)\bigr)$. As with compact quantum groups there are obvious examples provided by ordinary compact semigroups. Moreover, any compact quantum space admits a structure of a compact quantum semigroup, but there are interesting examples of quantum spaces (with natural quantum semigroup structures) which do not admit any compact quantum group structure.

\section{Quantum families of maps}

\subsection{Classical families of maps --- Jackson's theorem}\label{Jackson}

In classical topology the sets of maps between topological spaces are often themselves given a topology. A particular example is the \emph{compact-open} topology on the space of continuous maps from a topological space $X$ to a topological space $Y$ in which a basis of neighborhoods of a given $\phi\colon{X}\to{Y}$ is indexed by open sets $\mathcal{O}$ in $Y$ and compact sets $K$ in $X$ such that $\phi(K)\subset\mathcal{O}$ and consists of all other (continuous) $\psi\colon{X}\to{Y}$ which satisfy $\psi(K)\subset\mathcal{O}$. Thus $\C(X,Y)$ becomes a topological space and we can consider \emph{continuous families of maps $X\to{Y}$}, i.e.~continuous maps from a topological space $E$ to $\C(X,Y)$. The fundamental result about such families is Jackson's theorem:

\begin{theorem}[J.R.~Jackson]
Let $X$, $Y$ and $E$ be topological spaces such that $X$ is Hausdorff and $E$ is locally compact. For $\psi\in\C(X\times{E},Y)$ define $\sigma(\psi)$ as the mapping from $E$ to $\C(X,Y)$ given by
\[
\bigl(\sigma(\psi)(e)\bigr)(x)=\psi(x,e),\qquad(x\in{X},\;e\in{E}).
\]
Then $\sigma$ is a homeomorphism of $\C(X\times{E},Y)$ onto $\C\bigl(E,\C(X,Y)\bigr)$ with all spaces of maps topologized by their respective compact-open topologies.
\end{theorem}

In other words a continuous family of maps $X\to{Y}$ parametrized by points of $E$ is encoded in a single map $E\to\C(X,Y)$ and vice versa.

Before stating the definition of a quantum family of maps let us dwell for a little while on the \emph{universal property} of the space $\C(X,Y)$. Assume that $X$ is locally compact and let $\Gamma$ be the mapping $\C(X,Y)\times{X}\to{Y}$ defined by
\[
\Gamma(f,x)=f(x).
\]
It is known that $\Gamma$ is continuous. Moreover if $E$ is a locally compact space and $\psi\colon{X}\times{E}\to{Y}$ is continuous then there exists a unique continuous map $\lambda\colon{E}\to\C(X,Y)$ such that
\[
\Gamma\bigl(\lambda(e),x\bigr)=\psi(x,e).
\]
Indeed, the map $\lambda$ is simply $\sigma(\psi)$, where $\sigma$ is the homeomorphism from Jackson's theorem. The uniqueness follows from the fact that $\sigma$ is one-to-one. This universal property of the pair $\bigl(\C(X,Y),\Gamma\bigr)$ characterizes it uniquely. We will use this property to define the analog of $\C(X,Y)$ for quantum spaces.

\subsection{Definition of a quantum family of maps}

\begin{definition}\label{qfDef}
Let $\XX$, $\YY$ and $\EE$ be quantum spaces. A continuous \emph{quantum family of maps $\XX\to\YY$} parametrized by $\EE$ is an element
\[
\Psi\in\Mor\bigl(\C_0(\YY),\C_0(\XX)\tens\C_0(\EE)\bigr).
\]
\end{definition}

Let us analyze Definition \ref{qfDef} in more detail. The morphism $\Psi$ describes a continuous map from the quantum space corresponding to $\C(\XX)\tens\C(\EE)$ to $\YY$. In view of Jackson's theorem and the fact that tensor product corresponds to Cartesian product of spaces (cf.~Proposition \ref{cartprod}) we can interpret this as a family of maps $\XX\to\YY$ parametrized by the quantum space $\EE$. Specifying the definition of a quantum family of maps to the ``classical case'', i.e.~when all considered \cst-algebras are commutative, we obtain the standard notion of a continuous family of maps:

\begin{proposition}\label{classqfm}
Let $X,Y$ and $E$ be locally compact spaces and let $\Psi\in\Mor\bigl(\C_0(Y),\C_0(X)\tens\C_0(E)\bigr)$. Then there exists a unique family $\{\psi_e\}_{e\in{E}}$ of maps $X\to{Y}$ parametrized continuously by $E$ such that
\[
\Psi(f)(x,e)=f\bigl(\psi_e(x)\bigr)
\]
for all $x\in{X}$ and $e\in{E}$.
\end{proposition}

More generally, if only $\EE$ is a \emph{classical space} (which means that $\C_0(\EE)$ is commutative) then there exists a unique family $\{\Psi_e\}_{e\in{\EE}}$ of morphisms from $\C_0(\YY)$ to $\C_0(\XX)$ parametrized continuously (in a natural topology on the set of morphisms) by points of $\EE$.

Since a quantum family of maps is simply a morphism from a \cst-algebra to a tensor product of two \cst-algebras, it is not our aim to derive a general theory of such objects, as it would not likely give non-trivial results. We will rather focus on particular examples of quantum families, especially those possessing certain \emph{universal properties} (see Lecture II).

Before we deal with quantum families of maps with such additional properties let us define the composition of quantum families of maps. This is a direct generalization of the notion of composition of classical families of maps.

\begin{definition}
Let $\XX_1,\XX_2,\XX_3,\DD_1$ and $\DD_2$ be quantum spaces. Consider families of maps
\[
\begin{split}
\Psi_1&\in\Mor\bigl(\C_0(\XX_2),\C_0(\XX_1)\tens\C_0(\DD_1)\bigr),\\
\Psi_2&\in\Mor\bigl(\C_0(\XX_3),\C_0(\XX_2)\tens\C_0(\DD_2)\bigr).
\end{split}
\]
The \emph{composition} of $\Psi_1$ and $\Psi_2$ is the quantum family of maps
\[
\Psi_1\vt\Psi_2\in\Mor\bigl(\C_0(\XX_3),\C_0(X_1)\tens\C_0(\DD_1)\tens\C_0(\DD_2)\bigr)
\]
defined by
\[\Label{compo}
\Psi_1\vt\Psi_2=(\Psi_1\tens\id)\comp\Psi_2.
\]
\end{definition}

Note that the composition of morphisms in \eqref{compo} is defined via extending $(\Psi_1\tens\id)$ to a $*$-homomorphism from $\M\bigl(\C_0(\XX_2)\tens\C_0(\DD_2)\bigr)$ to $\M\bigl(\C_0(X_1)\tens\C_0(\DD_1)\tens\C_0(\DD_2)\bigr)$.

It is easy to see that if the considered quantum families are classical (i.e.~the algebras $\C(\DD_1)$ and $\C(\DD_2)$ are commutative) then composition of quantum families corresponds exactly to the family of all compositions of $*$-homomorphisms from the families $\Psi_1$ and $\Psi_2$. The operation of taking a composition of quantum families is associative in the following sense: let $\XX_1,\XX_2,\XX_3,\XX_4,\DD_1,\DD_2$ and $\DD_3$ be quantum spaces and let
\[
\begin{split}
\Psi_1&\in\Mor\bigl(\C_0(\XX_2),\C_0(\XX_1)\tens\C_0(\DD_1)\bigr),\\
\Psi_2&\in\Mor\bigl(\C_0(\XX_3),\C_0(\XX_2)\tens\C_0(\DD_2)\bigr),\\
\Psi_3&\in\Mor\bigl(\C_0(\XX_4),\C_0(\XX_3)\tens\C_0(\DD_3)\bigr).
\end{split}
\]
Then $(\Psi_1\vt\Psi_2)\vt\Psi_3=\Psi_1\vt(\Psi_2\vt\Psi_3)$.

From now on we will concentrate on \emph{compact} quantum spaces (the considered \cst-algebras will be unital). The reasons for this are the following:
\begin{itemize}
\item it is technically easier, in particular morphisms become simply unital $*$-homomorphisms,
\item in all considered examples the interesting \cst-algebras will be unital,
\item the universal properties formulated in the category of compact quantum spaces in fact hold in the category of all quantum spaces and the technical details are not at all difficult, while the notation is much heavier.
\end{itemize}
Therefore we will from now on even refrain from writing $\Psi\in\Mor(\sA,\sB)$ favoring the notation $\Phi\colon\sA\to\sB$.

\section{Quantum families of all maps, existence question}

\begin{definition}\label{DefAll}
Let $\XX$ and $\YY$, $\EE$ be quantum spaces and let $\Phi\colon\C(\YY)\to\C(\XX)\tens\C(\EE)$ be a quantum family of maps. We say that
\begin{itemize}
\item $\Phi$ is the \emph{quantum family of all maps} from $\XX$ to $\YY$ and
\item $\EE$ is the \emph{quantum space of all maps} from $\XX$ to $\YY$
\end{itemize}
if for any quantum space $\DD$ and any quantum family $\Psi\colon\C(\YY)\to\C(\XX)\tens\C(\DD)$ there exists a unique $\Lambda\colon\C(\EE)\to\C(\DD)$ such that the diagram
\[
\xymatrix{
\C(\YY)\ar[rr]^-{\Phi}\ar@{=}[d]&&\C(\XX)\tens\C(\EE)\ar[d]^{\id\tens\Lambda}\\
\C(\YY)\ar[rr]^-{\Psi}&&\C(\XX)\tens\C(\DD)}
\]
commutes.
\end{definition}

The motivation for Definition \ref{DefAll} comes from the universal property of the space $\C(X,Y)$ described at the end of Subsection \ref{Jackson}. Indeed, upon specializing to the classical situation ($\XX,\YY,\EE$ are taken to be classical spaces as well as all possible spaces $\DD$), Definition \ref{DefAll} turns out to be exactly the definition of the space $\C(\XX,\YY)$ with its compact open topology). This result is not entirely trivial.

The immediate question now is whether, given two quantum space $\XX$ and $\YY$ the quantum space of all maps from $\XX$ to $\YY$ (and the corresponding quantum family) exists. In the next theorem we use the notion of a \emph{finitely generated unital \cst-algebra}. By this we mean an algebra which is a quotient of the full group \cst-algebra of a free group on finite number of generators.

\begin{theorem}\label{SLW}
Let $\XX$ and $\YY$ be quantum spaces such that $\C(\XX)$ is finite dimensional and $\C(\YY)$ is finitely generated and unital. Then the quantum space of all maps $\XX\to\YY$ exists. Moreover this quantum space is compact.
\end{theorem}

One way to make sure the assumptions of Theorem \ref{SLW} are satisfied is to assume that both $\C(\XX)$ and $\C(\YY)$ are finite dimensional. Indeed many very interesting examples come from the situation when $\XX$ and $\YY$ are the same \emph{finite} quantum space $\MM$. Contrary to the classical situation, examples of such spaces can be very intriguing.

\begin{example}\label{Ex1}
Let $\MM$ be the classical two point space: $\C(\MM)=\CC^2$. Clearly the \emph{classical} space of all maps $\MM\to\MM$ has four points. However, the \emph{quantum} space of all maps $\MM\to\MM$ is much more interesting. It turns out to be a quantum space $\EE$ such that
\[
\C(\EE)=\bigl\{f\in\C\bigl([0,1],M_2(\CC)\bigr)\st{f(0),f(1)}\text{ are diagonal}\bigr\}.
\]
The quantum family of all maps $\MM\to\MM$ is $\Phi\colon\CC^2\to\CC^2\tens\C(\EE)$ given by
\[
\Phi\bigl(\left[\begin{smallmatrix}1\\0\end{smallmatrix}\right]\bigr)=
\left[\begin{smallmatrix}1\\0\end{smallmatrix}\right]\tens{P}+
\left[\begin{smallmatrix}0\\1\end{smallmatrix}\right]\tens{Q},
\]
where
\[
P(t)=\begin{bmatrix}0&0\\0&1\end{bmatrix},\quad
Q(t)=\tfrac{1}{2}\begin{bmatrix}
1-\cos{2\pi{t}}&\mathrm{i}\sin{2\pi{t}}\\
-\mathrm{i}\sin{2\pi{t}}&1+\cos{2\pi{t}}
\end{bmatrix},\qquad(t\in[0,1]).
\]
The infinite dimensionality of $\C(\EE)$ justifies saying that the quantum space of all maps from a two point space to itself is infinite.

The explanation of these formulas is quite simple: in order to determine a $*$-homomorphism from $\CC^2$ to $\CC^2\tens\sA$, where $\sA$ is some unital \cst-algebra, we only need to define its value on $\left[\begin{smallmatrix}1\\0\end{smallmatrix}\right]$. This value can be written as
\[
\left[\begin{smallmatrix}1\\0\end{smallmatrix}\right]\tens{a}+
\left[\begin{smallmatrix}0\\1\end{smallmatrix}\right]\tens{b},
\]
and $a,b\in\sA$ must be projections because $\left[\begin{smallmatrix}1\\0\end{smallmatrix}\right]$ is. The \cst-algebra $\C(\EE)$ described above is simply the universal \cst-algebra generated by two projections without any relations.
\end{example}

\section*{{\Large Lecture II}}
\setcounter{section}{0}
\renewcommand{\thesection}{II.\arabic{section}}

\begin{proof}[Idea of proof of Theorem \ref{SLW}]
Begin by proving the theorem for $\C(\YY)=\cst(\FF_n)$. If $\C(\XX)$ has the decomposition
\[
\C(\XX)=\bigoplus_{i=1}^NM_{n_i}(\CC)
\]
then we define $\EE_n$ so that $\C(\EE_n)$ is generated by
\[
\bigl\{u^{i,p}_{k,l}\st{1\leq{i}\leq{N}},\;{1\leq{k,l}\leq{n_i},\;{1\leq{p}\leq{n}}}\bigr\}
\]
with relations that the matrices
\[
u^{i,p}=\begin{bmatrix}u^{i,p}_{1,1}&\dotsm&u^{i,p}_{1,m_i}\\
\vdots&\ddots&\vdots\\u^{i,p}_{m_i,1}&\dotsm&u^{i,p}_{m_i,m_i}\end{bmatrix},\qquad(1\leq{i}\leq{N},\:1\leq{p}\leq{n})
\]
be unitary. The we define a map $\Phi_n$ as then the unique $*$-homomorphism $\cst(\FF_n)\to\C(\XX)\tens\C(\EE_n)$ by
\[
\Phi(v_p)=(u^{1,p},\dotsc,u^{N,p})\in\bigoplus_{i=1}^NM_{m_i}\bigl(\C(\EE_n)\bigr)\cong\C(\XX)\tens\C(\EE_n),
\]
where $\{v_1,\dotsc,v_n\}$ are the generators of $\FF_n$. The pair $\bigl(\C(\EE_n),\Phi_n)$ has the universal property that for any quantum space $\DD$ and any quantum family $\Psi\colon\cst(\FF_n)\to\C(\XX)\tens\C(\DD)$ there exists a unique $\Lambda_n\colon\C(\EE_n)\to\C(\DD)$ such that the diagram
\[
\xymatrix{
\mathrm{C}^*(\FF_n)\ar[rr]^-{\Phi_n}\ar@{=}[d]&&\C(\XX)\tens\C(\EE_n)\ar[d]^{\id\tens\Lambda_n}\\
\mathrm{C}^*(\FF_n)\ar[rr]^-{\Psi}&&\C(\XX)\tens\C(\DD)}
\]
commutes. Now we take $\YY$ such that $\C(\YY)$ is a quotient of $\cst(\FF_n)$ and let $\pi\colon\cst(\FF_n)\to\C(\YY)$ be the quotient map. Define $\EE$ so that $\C(\EE)$ is the quotient of $\CC(\EE_n)$ by the ideal generated by the set
\[
\bigl\{(\omega\tens\id)\Phi_n(x)\st{x\in\ker\pi},\;\omega\in{\C(\XX)^*}\bigr\}
\]
and let $\lambda\colon\C(\EE_n)\to\C(\EE)$ be the quotient map. Now it can be shown that there exists a map $\Phi\colon\C(\YY)\to\C(\XX)\tens\C(\EE)$ such that diagram
\begin{equation}\label{bPhi}
\xymatrix{
{\cst}(\FF_n)\ar[rr]^-{\Phi_n}\ar[d]_{\pi}&&\C(\XX)\tens\C(\EE_n)\ar[d]^{\id\tens\lambda}\\
\C(\YY)\ar[rr]^-{\Phi}&&\C(\XX)\tens\C(\EE)}
\end{equation}
is commutative.

The pair $\bigl(\C(\EE),\Phi\bigr)$ will be the quantum space and family of all maps $\XX\to\YY$. The universal property follows from the fact that if $\DD$ is a quantum space and $\Psi\colon\C(\YY)\to\C(\XX)\tens\C(\DD)$ is a quantum family of maps then by the first step of the proof, there exists a unique $\Lambda_n\colon\C(\EE_n)\to\C(\DD)$ such that
\[
\Psi\comp\pi=(\id\tens\Lambda_n)\comp\Phi_n.
\]
Then it can be shown that there exists a unique $\Lambda\colon\C(\EE)\to\C(\DD)$ such that $\Lambda_n=\Lambda\comp\lambda$ (this follows form the construction of $\C(\EE)$ and $\lambda$). Therefore the diagram
\[
\xymatrix
{
\C(\XX)\tens\C(\EE_n)\ar[rr]^-{\id\tens\Lambda_n}\ar[rd]_-{\id\tens\lambda}&&\C(\XX)\tens\C(\DD)\\
&\C(\XX)\tens\C(\EE)\ar[ur]_-{\id\tens\Lambda}
}
\]
commutes. It follows that the diagram
\begin{equation}\label{tego}
\xymatrix{{\cst(\FF_n)}\ar[rr]^-{\Phi_n}\ar[dd]_{\pi}&&\C(\XX)\tens\C(\EE_n)\ar[dl]_{\id\tens\lambda}
\ar[dd]^{\id\tens\Lambda_n}\\
&\C(\XX)\tens\C(\EE)\ar[rd]^-{\id\tens\Lambda}\\
\C(\YY)\ar[rr]_-{\Psi}\ar[ur]^-{\Phi}&&\C(\XX)\tens\C(\DD)}
\end{equation}
commutes. The diagram
\[
\xymatrix{\C(\YY)\ar[rr]^-{\Phi}\ar@{=}[d]&&\C(\XX)\tens\C(\EE)\ar[d]^{\id\tens\Lambda}\\\C(\YY)\ar[rr]^-{\Psi}&&\C(\XX)\tens\C(\DD)}
\]
is simply a part of \eqref{tego}.
\end{proof}

The fact that the quantum space of all maps $\XX\to\YY$ might not exists should rather be interpreted as our inability to see it. More precisely, the approach to non-commutative topology we have adopted (via \cst-algebras) should be viewed more as ``non-commutative topology of \emph{locally compact} quantum spaces''. It is clear that the space of continuous mappings with compact open topology often fails to be locally compact. This happens also in the classical case.

\section{Quantum semigroup structure on the quantum space of all maps}

Choose a finite quantum space $\MM$. Then the quantum space of all maps $\MM\to\MM$ exists. One expects that this space should be a semigroup under composition.

\begin{theorem}
Let $\MM$ be a finite quantum space and let $\EE$ be the quantum space of all maps $\MM\to\MM$. Let
\[
\Phi\colon\C(\MM)\longrightarrow\C(\MM)\tens\C(\EE)
\]
be the quantum family of all maps $\MM\to\MM$. Then there exists a unique
\[
\Delta\colon\C(\EE)\longrightarrow\C(\EE)\tens\C(\EE)
\]
such that $\Phi\vt\Phi=(\id\tens\Delta)\comp\Phi$ i.e.~the diagram
\[
\xymatrix{
\C(\MM)\ar[rr]^-{\Phi}\ar@{=}[d]&&\C(\MM)\tens\C(\EE)\ar[d]^{\id\tens\Delta}\\
\C(\MM)\ar[rr]^-{\Phi\vt\Phi}&&\C(\MM)\tens\C(\EE)\tens\C(\EE)}
\]
is commutative. Moreover the $*$-homomorphism $\Delta$ is coassociative:
\[
(\Delta\tens\id)\comp\Delta=(\id\tens\Delta)\comp\Delta
\]
and there exists a unique character $\epsilon$ of $\C(\EE)$ such that
\[
(\id\tens\epsilon)\comp\Phi=\id
\]
which also satisfies $(\epsilon\tens\id)\comp\Delta=(\id\tens\epsilon)\comp\Delta=\id$.
\end{theorem}

In other words the quantum space of all maps $\MM\to\MM$ has a natural structure of a compact quantum semigroup. The character $\epsilon$ corresponds to the unit element of this semigroup and is called a \emph{counit}.

\begin{example}
Let, as in Example \ref{Ex1}, $\C(\MM)=\CC^2$ and let $\EE$ be the quantum space of all maps $\MM\to\MM$. We know that $\C(\EE)$ is generated by the two distinguished elements $P$ and $Q$ (see Example \ref{Ex1}). It is easy to find that the quantum semigroup structure on $\EE$ is given by the comultiplication:
\[
\Delta(P)=P\tens{P}+(\I-P)\tens{Q},\qquad\Delta(Q)=Q\tens{P}+(\I-Q)\tens{Q}.
\]
Also $\epsilon(P)=1$ and $\epsilon(Q)=0$.
\end{example}

It is natural to expect that the quantum semigroup of all maps $\MM\to\MM$ is not a quantum \emph{group}. Indeed we have

\begin{proposition}
Unless $\MM$ is a (classical) one point space, the quantum semigroup of all maps $\MM\to\MM$ is not a quantum group.
\end{proposition}

\section{Subsemigroups: quantum families preserving a state, quantum commutants}

\subsection{Quantum semigroups preserving a state}\label{state}

If $M$ is a classical topological space and $\mu$ a measure on $M$ then the set of those continuous maps $M\to{M}$ that preserve $\mu$ is a subsemigroup of $\C(M,M)$ (this fact has little to do with topology). We are now looking for an analog of this phenomenon in the non-commutative setting.

Let $\XX$ be a quantum space and let $\omega$ be a state on $\C(\XX)$ (i.e.~a positive linear functional of norm 1 --- an analog of a probability measure on $\MM$). We say that a quantum family of maps $\Psi\colon\C(\XX)\to\C(\XX)\tens\C(\DD)$ \emph{preserves} $\omega$ if
\[
(\omega\tens\id)\Psi(a)=\omega(a)\I
\]
for all $a\in\C(\XX)$. This is a direct generalization of the standard notion of a family of maps preserving a given (probability) measure.

\begin{theorem}\label{inv_omega}
Let $\MM$ be a finite quantum space and $\omega$ a state on $\C(\MM)$. Then
\begin{enumerate}
\item
there exists a unique quantum family
\[
\Phi_\omega\colon\C(\MM)\longrightarrow\C(\MM)\tens\C(\WW)
\]
such that for any quantum family
\[
\Psi\colon\C(\MM)\longrightarrow\C(\MM)\tens\C(\DD)
\]
preserving $\omega$ there exists a unique $\Theta\colon\C(\WW)\to\C(\DD)$ such that
\[
\xymatrix{
\C(\MM)\ar[rr]^-{\Phi_\omega}\ar@{=}[d]&&\C(\MM)\tens\C(\WW)\ar[d]^{\id\tens\Theta}\\
\C(\MM)\ar[rr]^-{\Psi}&&\C(\MM)\tens\C(\DD)}
\]
\item $\Phi_\omega$ preserves $\omega$,
\item $\WW$ has a structure of a compact quantum semigroup with unit element given by comultiplication $\Delta_\omega\colon\C(\WW)\to\C(\WW)\tens\C(\WW)$ and if $\EE$ is the quantum semigroup of all maps $\MM\to\MM$ then the canonical map $\Lambda\colon\C(\EE)\to\C(\WW)$ intertwines the comultiplications of $\C(\EE)$ and $\C(\WW)$:
\[
\xymatrix{
\C(\EE)\ar[rr]^-{\Delta}\ar[d]_{\Lambda}&&\C(\EE)\tens\C(\EE)\ar[d]^{\Lambda\tens\Lambda}\\
\C(\WW)\ar[rr]_-{\Delta_\omega}&&\C(\WW)\tens\C(\WW)
}
\]
\end{enumerate}
\end{theorem}

Theorem \ref{inv_omega} says that there exists the \emph{quantum family of all maps $\MM\to\MM$ preserving $\omega$}. The quantum space $\WW$ parameterizing this family is a compact quantum semigroup, in fact, a subsemigroup of the quantum semigroup $\EE$ of all maps $\MM\to\MM$. The \cst-algebra $\CC(\WW)$ is constructed as the quotient of $\CC(\EE)$ by the ideal generated by
\[
\bigl\{(\omega\tens\id)\Phi(m)-\omega(m)\I\st{m}\in\C(\MM)\bigr\}.
\]

\begin{example}\label{Ex2}
Let $\C(\MM)=M_2(\CC)$ and choose a parameter $q\in]0,1]$. Let $\omega_q$ be the state on $\C(\MM)$ defined by
\[
\omega_q\left(\left[\begin{smallmatrix}a&b\\c&d\end{smallmatrix}\right]\right)=\tfrac{a+q^2d}{1+q^2}
\]
The quantum semigroup $\WW$ of all maps $\MM\to\MM$ preserving $\omega_q$ can be described as follows: $\C(\WW)$ is generated by three elements $\beta$, $\gamma$ and $\delta$ such that
\[
\begin{aligned}
q^4\delta^*\delta+\gamma^*\gamma+q^4\delta\delta^*+\beta\beta^*&=\I,
&\beta\gamma&=-q^4\delta^2,\\
\beta^*\beta+\delta^*\delta+\gamma\gamma^*+\delta\delta^*&=\I,
&\gamma\beta&=-\delta^2,\\
\gamma^*\delta-q^2\delta^*\beta+\beta\delta^*-q^2\delta\gamma^*&=0,
&\beta\delta&=q^2\delta\beta,\\
q^4\delta\delta^*+\beta\beta^*+q^2\gamma\gamma^*+q^2\delta\delta^*&=\I,&\delta\gamma&=q^2\gamma\delta,\\
q^4\delta^*\delta+\gamma^*\gamma+q^2\beta^*\beta+q^2\delta^*\delta&=q^2\I.
\end{aligned}
\]
The comultiplication $\Delta_{\omega_q}\colon\C(\WW)\to\C(\WW)\tens\C(\WW)$ is given by
\[
\begin{split}
\Delta_{\omega_q}(\beta)&=q^4\delta\gamma^*\tens\delta-q^2\beta\delta^*\tens\delta+\beta\tens\beta
+\gamma^*\tens\gamma-q^2\delta^*\beta\tens\delta+\gamma^*\delta\tens\delta,\\
\Delta_{\omega_q}(\gamma)&=q^4\gamma\delta^*\tens\delta-q^2\delta\beta^*\tens\delta+\gamma\tens\beta
+\beta^*\tens\gamma-q^2\beta^*\delta\tens\delta+\delta^*\gamma\tens\delta,\\
\Delta_{\omega_q}(\delta)&=-q^2\gamma^*\gamma\tens\delta-q^2\delta\delta^*\tens\delta+\delta\tens\beta
+\delta^*\tens\gamma+\beta^*\beta\tens\delta+\delta^*\delta\tens\delta.
\end{split}
\]
The counit $\epsilon$ maps $\gamma$ and $\delta$ to $0$ and $\beta$ to $1$.

The proof of the above facts is mainly computational. One first observes that $M_2(\CC)$ is the universal \cst-algebra generated by one element $n$ such that
\[
n^2=0,\qquad{nn^*+n^*n=\I}.
\]
(in fact we can specify $n=\bigl[\begin{smallmatrix}0&1\\0&0\end{smallmatrix}\bigr]$). Then the mapping $\Phi_{\omega_q}$ will map $n$ to a matrix
\[
\begin{bmatrix}\alpha&\beta\\\gamma&\delta\end{bmatrix}
\]
of elements of $\C(\WW)$. The relations between matrix elements are derived from
\[
\Phi_{\omega_q}(n)^2=0,\qquad\Phi_{\omega_q}(n)\Phi_{\omega_q}(n)^*+\Phi_{\omega_q}(n)^*\Phi_{\omega_q}(n)=\I
\]
and
\[
\begin{split}
(\omega_q\tens\id)\Phi_{\omega_q}(n)&=0,\\
(\omega_q\tens\id)\Phi_{\omega_q}(n)&=0,\\
(\omega_q\tens\id)\bigl(\Phi_{\omega_q}(n)\Phi_{\omega_q}(n)^*\bigr)&=\tfrac{1}{1+q^2}\I,\\
(\omega_q\tens\id)\bigl(\Phi_{\omega_q}(n)^*\Phi_{\omega_q}(n)\bigr)&=\tfrac{q^2}{1+q^2}\I.
\end{split}
\]
In particular it follows from the relations that $\alpha=-q^2\delta$, so
\[
\Phi_{\omega_q}(n)=\begin{bmatrix}-q^2\delta&\beta\\\gamma&\delta\end{bmatrix}
\in{M_2(\CC)}\tens\C(\WW)=\C(\MM)\tens\C(\WW).
\]
\end{example}

\subsection{Quantum commutants}

In Subsection \ref{state} we found that in the non-commutative setting there exists an analog of the space of all continuous maps of a given space $M$ into itself preserving a fixed probability measure and it is in an appropriate sense a subsemigroup of the set of all maps $M\to{M}$. Another similar object is the family of all maps \emph{commuting} with a fixed family of maps.

\begin{definition}
Let $\MM$ be a finite quantum space and let
\[
\Psi_1\colon\C(\MM)\to\C(\MM)\tens\C(\DD_1),\qquad
\Psi_2\colon\C(\MM)\to\C(\MM)\tens\C(\DD_2)
\]
be two quantum families of maps. We say that $\Psi_1$ and $\Psi_2$ \emph{commute} if
\[
(\id\tens\sigma)\comp(\Psi_1\vt\Psi_2)=\Psi_2\vt\Psi_1,
\]
where $\sigma$ is the flip
\[
\C(\DD_1)\tens\C(\DD_2)\ni{x\tens{y}}\longmapsto{y\tens{x}}\in\C(\DD_2)\tens\C(\DD_1).
\]
\end{definition}

\begin{theorem}
Let $\MM$ be a finite quantum space and $\Psi\colon\C(\MM)\to\C(\MM)\tens\C(\DD)$ a quantum family of maps $\MM\to\MM$. Then
\begin{enumerate}
\item there exists a unique quantum family
\[
\Phi_\Psi\colon\C(\MM)\longrightarrow\C(\MM)\tens\C(\UU)
\]
such that for any quantum family $\Xi\colon\C(\MM)\to\C(\MM)\tens\C(\PP)$ commuting with $\Psi$ there exists a unique $\Theta\colon\C(\UU)\to\C(\PP)$ making the diagram
\[
\xymatrix{
\C(\MM)\ar[rr]^-{\Phi_\Psi}\ar@{=}[d]&&\C(\MM)\tens\C(\UU)\ar[d]^{\id\tens\Theta}\\
\C(\MM)\ar[rr]^-{\Xi}&&\C(\MM)\tens\C(\PP)}
\]
commutative.
\item $\Phi_\Psi$ commutes with $\Psi$,
\item $\UU$ is a compact quantum semigroup with unit and if $\EE$ is the quantum semigroup of all maps $\MM\to\MM$ then the canonical map $\Lambda\colon\C(\EE)\to\C(\UU)$ intertwines the comultiplications of $\C(\EE)$ and $\C(\UU)$.
\end{enumerate}
\end{theorem}

Of course any classical family of morphisms $\C(\MM)\to\C(\MM)$ can be interpreted as a quantum family (parametrized by a classical space) we can consider quantum commutants also of classical families --- even consisting of a single morphism.

\begin{example}
Let, as in Example \ref{Ex2}, $\MM$ be a quantum space such that $\C(\MM)=M_2(\CC)$. Let $\UU$ be the commutant of the (classical) family of maps $\MM\to\MM$ consisting of the single automorphism of $\C(\MM)$ defined by
\[
\psi\colon\left[\begin{smallmatrix}a&b\\c&d\end{smallmatrix}\right]\longmapsto
\left[\begin{smallmatrix}d&c\\b&a\end{smallmatrix}\right].
\]
This family is described in our language by the $*$-homomorphism
\[
\Psi\colon\C(\MM)\longrightarrow\C(\MM)\tens\CC
\]
given by $\Psi(m)=\psi(m)\tens{1}$. Let $\UU$ be the quantum commutant of $\Psi$. The \cst-algebra $\C(\UU)$ is generated by $\alpha,\beta$ and $\gamma$ with
\[
\beta=\beta^*,\quad\gamma=\gamma^*
\]
and
\[
\begin{aligned}
\alpha^*\alpha+\gamma^2+\alpha\alpha^*+\beta^2&=\I,
&\alpha^2+\beta\gamma&=0,\\
\alpha^*\beta+\gamma\alpha^*+\alpha\gamma+\beta\alpha&=0,
&\alpha\beta+\beta\alpha^*&=0,\\
\gamma\alpha+\alpha^*\gamma&=0.
\end{aligned}
\]
The comultiplication $\Delta_\Psi\colon\C(\UU)\to\C(\UU)\tens\C(\UU)$ acts on generators as
\[
\begin{split}
\Delta_\Psi(\alpha)&=\I\tens\alpha+(\alpha^*\alpha+\gamma^2)\tens(\alpha^*-\alpha)+
\alpha\tens\beta+\alpha^*\tens\gamma,\\
\Delta_\Psi(\beta)&=(\alpha\gamma+\beta\alpha)\tens(\alpha-\alpha^*)+\beta\tens\beta+\gamma\tens\gamma,\\
\Delta_\Psi(\gamma)&=(\beta\alpha+\alpha\gamma)\tens(\alpha^*-\alpha)+\gamma\tens\beta+\beta\tens\gamma,
\end{split}
\]
It can be shown that $\UU$ is not a compact quantum group (with the comultiplication $\Delta_\Psi$).
\end{example}

\section{Further topics: quantum space of maps from a finite quantum space to a quantum semigroup and its structure}

So far we have concentrated on subsemigroups of the quantum semigroup of all maps from a given finite quantum space to itself. In classical topology there is also another way to form semigroups of maps. Namely one considers the space of all continuous maps form a fixed topological space to a topological semigroup. This time the group operation is not composition of maps, but a pointwise multiplication.

A non-commutative analog of this construction was investigated first by M.M.~Sadr. His idea is contained in the following theorem:

\begin{theorem}\label{qms}
Let $K$ be a finite set and let $\sS$ be a quantum semigroup with comultiplication $\Delta_\sS\colon\C(\sS)\to\C(\sS)\tens\C(\sS)$. Assume that $\C(\sS)$ is finitely generated and unital. Then let $\HH$ be the quantum space of all maps $K\to\sS$ and let
\[
\Phi\colon\C(\sS)\to\C(K)\tens\C(\HH)
\]
be the quantum family of all these maps. Then $\HH$ admits a structure of a compact quantum semigroup described by $\Delta_\HH\colon\C(\HH)\to\C(\HH)\tens\C(\HH)$ defined by the diagram
\[
\xymatrix
{
\C(\sS)\ar[rr]^-\Phi\ar[d]_{\Delta_\sS}&&\C(K)\tens\C(\HH)\ar[d]^{\id\tens\Delta_\HH}\\
\C(\sS)\tens{\C(\sS)}\ar[d]_{\Phi\tens\Phi}&&\C(K)\tens\C(\HH)\tens\C(\HH)\\
\C(K)\tens\C(\HH)\tens\C(K)\tens\C(\HH)\ar[rr]_{\id\tens\sigma\tens\id}
&&\C(K)\tens\C(K)\tens\C(\HH)\tens\C(K)\ar[u]_{\mu\tens\id\tens\id}
}
\]
where $\sigma$ is the flip $\C(\HH)\tens\C(K)\to\C(K)\tens\C(\HH)$ and $\mu\colon\C(K)\tens\C(K)\to\C(K)$ is the multiplication map.
\end{theorem}

The crucial remarks are:

\begin{itemize}
\item Upon substituting a classical semigroup $S$ for $\sS$ and taking the \emph{classical} space of all maps $K\to{S}$ instead of $\HH$, we obtain exactly the structure of the semigroup of all maps $K\to{S}$.
\item As stated, Theorem \ref{qms} can only be considered for \emph{classical} spaces $K$. Indeed, if we take a non-commutative \cst-algebra instead of $\C(K)$ the multiplication map $\mu$ is no longer an algebra homomorphism.
\end{itemize}

With a more detailed analysis of the structure of $\bigl(\C(\HH),\Delta_\HH\bigr)$ shows the following:

\begin{proposition}\label{iff}
Let $K$ be a finite set and let $\sS$ be a quantum semigroup as in Theorem \ref{qms}. Let $\HH$ be the quantum space of all maps ${K}\to\sS$ with the quantum semigroup structure as described in Theorem \ref{qms}. Then $\HH$ is a quantum group if and only if $\sS$ is.
\end{proposition}

This fact is quite puzzling in view of the following example:

\begin{example}
Let $\KK$ be the quantum space defined by $\C(\KK)=M_2(\CC)$ and let $\sS$ be the classical group $\ZZ_2$. Let $\HH$ be the quantum space of all maps $\KK\to\sS$. Then $\C(\HH)$ is isomorphic to the universal unital \cst-algebra generated by three elements $p,q$ and $z$ with the relations
\begin{align*}
&&&&p&=p^*,&p&=p^2+z^*z,&zp&=(\I-q)z,&&&&\\
&&&&q&=q^*,&q&=q^2+zz^*.
\end{align*}
It can be shown that this \cst-algebra does not admit a compact quantum group structure. Therefore one cannot hope for an analog of Proposition \ref{iff} when $K$ is not a classical space (since $\ZZ_2$ is a group, hence a quantum group, and $\HH$ does not admit a quantum group structure).
\end{example}


\begin{thebibliography}{6}
\bibitem{qs}
{\sc P.M.~So{\l}tan}: Quantum families of maps and quantum semigroups on finite quantum spaces. \emph{J.~Geom.~Phys.} \textbf{59} (2009) 354--368.
\bibitem{kom}
{\sc P.M.~So{\l}tan}: Examples of quantum commutants. In "Interactions of Algebraic \& Coalgebraic Structures (Theory and Applications)", \emph{Arab J.~Sci.~Eng.} \textbf{33} no.~2C (2008), 447--457.
\bibitem{apqs}
{\sc P.M.~So{\l}tan}: On quantum semigroup actions on finite quantum spaces. \emph{Infin.~Dimens.~Anal.~Qu.} \textbf{12} no.~3 (2009), 503--509.
\bibitem{qmqs}
{\sc P.M.~So{\l}tan}: On quantum maps into quantum semigroups. To appear in \emph{Houston Journal of Mathematics}. Available online at \texttt{ 	arXiv:1010.3379 [math.OA]}.
\end{thebibliography}
\end{document}